\documentclass{article}

\usepackage{cite}

\usepackage{amsmath,amsfonts,amsthm,amssymb,amscd,dsfont}

\numberwithin{equation}{section}

\usepackage[all,knot]{xy}
\xyoption{arc}
\newtheorem{theorem}{Theorem}[section]
\newtheorem{lemma}{Lemma}[section]

\newtheorem{corollary}[theorem]{Corollary}

\theoremstyle{definition}

\newtheorem{problem}[]{Problem}
\newtheorem{proof*}{Proof}

\theoremstyle{remark}
\newtheorem{remark}{Remark}[section]

\newcommand {\RR} {\mathbb R}
\newcommand {\NN} {\mathbb N}

\begin{document}

\title{On one extremal property of a regular simplex}

\author{Vladislav Babenko, Yuliya Babenko, \\ Nataliya Parfinovych, and Dmytro Skorokhodov}

\date{}

\maketitle

{ \centerline{\bf Abstract}

In this paper, we show that the $L_p$-error of asymmetric linear
approximation of the quadratic function $Q({\mathbf
x})=\sum_{j=1}^{d}x_j^2$ on simplices in $\RR^d$ of fixed volume
is minimized on regular simplices.}

\section{Introduction}\label{sec1}

One important problem in geometry is to study the approximation
(in a specified metric) of smooth convex bodies by various
polytopes. For instance, questions of approximating convex bodies
by inscribed or circumscribed polytopes, by polytopes with fixed
number of vertices or faces, by polytopes of best approximation,
etc. have been studied in this direction. After some occasional
results in the plane, the book of Fejes~Toth~\cite{Toth} was the
first to provide a large number of problems, ideas and results on
polytopal approximation in dimensions two and three, concentrating
specifically on extremal properties of regular polytopes. Many
extensions have been made afterwards to higher dimensions, other
metrics, etc. (see~\cite{boro,Boro1,Gruber,Handbook} and references
therein).

On the other hand, the question on approximation of functions,
defined on a polytope, by piecewise linear functions, generated
with the help of triangulations, in $L_p$-metrics is of great
importance in approximation theory. The problems of interpolation,
best and best one-sided approximation of functions by linear
splines (i.e., piecewise linear functions) have been considered.
The question of finding optimal adaptive triangulations, i.e.,
triangulations that depend on the function being approximated
(see, for instance,~\cite{Nira_Dyn_1,Nira_Dyn_2}), is of
particular interest for applications. In order to construct
adaptive triangulations many authors took into account the Hessian
of the function (or curvature of its graph)
(see~\cite{BBLS,BBS,chen_1,Daz3,huang,Nadler_86}).

Note that the construction of the best (in a specified sense)
polytope for an individual convex body, or construction of the
best triangulation for approximation of a specific individual
function, is an extremely difficult problem. The above-mentioned
papers have dealt only with asymptotically optimal sequences of
polytopes or asymptotically optimal sequences of triangulations.

One possible method to construct an asymptotically optimal
sequence of polytopes or triangulations is as follows. As the
first step, we construct an intermediate approximation of the
convex body surface (or the function, respectively) by piecewise
quadratic surface (function), and then solve the problem of
approximating the piecewise quadratic surface (function) by
piecewise linear ones. The latter, in turn (at least for $d=2$),
is equivalent to solving the following optimization problem (we
shall give its statement for approximation of functions).

Let a quadratic function $Q$ and a simplex $\mathcal{T}$
in $\RR^d$ of unit volume be given. We shall consider the best
$L_p$-approximation of function $Q$ by linear functions defined on
$\mathcal{T}$ (or the best one-sided approximation, which coincide
with the deviation of interpolant for positively definite $Q$),
and the problem is to find a simplex $\mathcal{T}^*$, for which
the corresponding error is minimal. (The known solutions of this
problem are listed in Section~3.)

Therefore, in a number of questions of geometry and approximation
theory, it is important to find a simplex of fixed (unit) volume
such that the error of the best approximation of a given quadratic
function on this simplex in a specified metric or the best
approximation with constraints (for instance, one-sided) is
minimized.

In approximation theory, there exists a tool to view both the
problem of finding the best approximation without constraints and
the problem of finding the best approximation with constraints
``under one umbrella''. The latter can be viewed as the best
approximation in the spaces with asymmetric norm or so-called
$(\alpha,\beta)$-approximation (see, for
example,~\cite{Bab_82,Bab_83,Korn}), when positive and negative
parts of the difference between function and the approximant are
``weighted'' differently. Such type of approximations are of a
separate interest, since they can be considered as the problems of
approximation with non-strict constraints (see below for more
precise statements), when constraints are allowed to be violated,
but the penalty for the violation is introduced into the error
measure. We think that such a method could also be interesting for
certain geometric problems.

Therefore, the purpose of this paper is to prove the optimality of
a regular simplex in the problem of minimizing (over the simplices
unit volume)\break the best ($\alpha,\beta$)-approximation in
$L_p$-metric of quadratic function\break $\sum_{j=1}^{d} x_j^2$ by
linear functions. Note that with the help of linear
transformations the solution of this problem allows us to obtain
the solution of analogous optimization problems for an arbitrary
positive definite quadratic\break form.

The paper is organized as follows. Section~2 contains definitions,
notations and rigorous statements of the problem we study and
closely related problems. Previous results and the main result of
the paper are given in Section~3. Section~4 is devoted to the
proof of the main result.

\section{Notations, definitions and statements of the problems}\label{sec2}

Let $d\in\NN$ and let $\RR^d$ be the space of points ${\mathbf
x}=(x_1,\ldots,x_d)$. Every point ${\mathbf x}\in\RR^d$ determines
(and is determined by) the row-vector with coordinates
$(x_1,\ldots,x_d)$, and we shall reserve the notation ${\mathbf
x}$ for such a vector. The Euclidean distance between points
${\mathbf a},{\mathbf b}\in\RR^d$ is defined as usually by
\[
\|{\mathbf a}-{\mathbf b}\|_{2}:=\left(\sum_{j=1}^d
(a_j-b_j)^2\right)^{1/2}.
\]
For a row-vector ${\mathbf x}$, let ${\mathbf x}^{\rm t}$ be the
column-vector transponent to ${\mathbf x}$. For a square matrix
${\mathbf J}$, we denote its transponent matrix by ${\mathbf
J}^{\rm t}$.

For a measurable bounded set $G\subset \RR^d$, let $L_{p}(G)$,
$1\le p\le \infty$, be the space of measurable and integrable in
the power $p$ (essentially bounded if $p=\infty$) and functions
$f:G\to\RR$ endowed with the usual norm
\[
\|f\|_{L_p(G)} : =\left\{\begin{array}{@{}ll@{}} \left(\int_{G} |f
({\mathbf x})|^p\,d{\mathbf
x}\right)^{\frac 1p} &\text{ if }1\le p < \infty, \\
{\rm esssup} \{|f ({\mathbf x}) |:{\mathbf x} \in G \} &\text{ if
} p =\infty.
\end{array}\right.
\]

Let $f\in L_p(G)$ and let a locally compact subset $H\subset
L_p(G)$ be given. Set $E(f;H)_{L_p(G)}$ to be the best
approximation of the function $f$ by $H$ in the $L_p$-norm, i.e.,
\[
E(f;H)_{L_p(G)}:=\inf\{\|f-u\|_{L_p(G)}\;:\;u\in H\}.
\]
In addition, set
\begin{equation}\label{trak}
E^{\pm}(f;H)_{L_p(G)}:=\inf\{\|f-u\|_{L_p(G)}\;:\;\pm u\le\pm f,\;u\in H\}.
\end{equation}
Quantity~(\ref{trak}) is called the best approximation from above
($E^-(f;H)_{L_p(G)}$) or below ($E^+(f;H)_{L_p(G)}$) of the
function $f$ by the subset $H$ in the $L_p$-norm. The quantities
$E^-(f;H)_{L_p(G)}$ and $E^+(f;H)_{L_p(G)}$ are called the best
one-sided approximations. For $\alpha,\beta>0$ and $f\in L_p(G)$,
let
\[
|f({\mathbf x})|_{\alpha,\beta}:=\alpha f_+({\mathbf x})+\beta f_-({\mathbf x}),
\]
where $g_{\pm}({\mathbf x}):=\max\{\pm g({\mathbf x});0\}$. Define
the asymmetric $L_p$-norm as follows:
\begin{align*}
\|f\|_{L_{p;\alpha,\beta}(G)} & :=  \|\alpha f_++\beta f_-\|_{L_p(G)} \\
& =  \begin{cases} \left(\int_{G} |f({\mathbf
x})|_{\alpha,\beta}^p\,d{\mathbf x}\right)^{\frac 1p} &\text{ if }
1\le p < \infty, \\[6pt]
{\rm esssup} \{|f({\mathbf x})|_{\alpha,\beta}:{\mathbf x} \in G
\} &\text{ if }p =\infty.
\end{cases}
\end{align*}
Asymmetric norms in connection with various problems in
approximation theory were considered in papers~\cite{Bab_82,Bab_87,Krein}
and books~\cite{KLD,Korn,KreinNudelman}. By\break $E(f;H)_{p;\alpha,\beta}$
denote the best ($\alpha,\beta$)-approximation~\cite{Bab_82,Korn}
of the function $f$ by $H$ in the $L_p$-norm, i.e.,
\[
E(f;H)_{L_{p;\alpha,\beta}(G)}=\inf\{\|f-u\|_{L_{p;\alpha,\beta}(G)}:\,u\in H\}.
\]
Note that for $\alpha \,{=}\, \beta \,{=}\, 1$, we have
$E(f;H)_{L_{p;1,1}(G)} \,{=}\, E(f;H)_{L_p(G)}$. Babenko~\cite{Bab_82}\break
proved that the following limit relations hold (see also
\cite[Theorem~1.4.10]{Korn}):
\begin{equation}\label{ter}
\begin{array}{l}
\displaystyle\lim_{\beta\to+\infty}E(f;H)_{L_{p;1,\beta}(G)}=E^+(f;H)_{L_p(G)}, \\
\displaystyle\lim_{\alpha\to+\infty}E(f;H)_{L_{p;\alpha,1}(G)}=E^-(f;H)_{L_p(G)}.
\end{array}
\end{equation}
This allows us to include the problem of the best unconstrained
approximation and the problem of the best one-sided approximation
into the family of problems of the same type, and consider them
from a general point of view (for more on this motivation,
see~\cite{Bab_83,Bab_84}). Because of the relation
\[
\|f-u\|_{p;1,\beta}^p=\|f-u\|_{p}^p + (\beta^p-1)\|(f-u)_-\|_p^p,
\quad \beta>1,
\]
the problem of the best $(1,\beta)$-approximation can be
considered as the problem of the best approximation with
non-strict constraint $f\le u$. This constraint is allowed to be
violated, but the penalty
\[
(\beta^p-1)\|(f-u)_-\|_{p}^p
\]
for its violation is introduced into the error measure. In what
follows, we shall allow the value $+\infty$ for $\alpha$ or
$\beta$, in that case identifying $E(f;H)_{L_{p;\alpha,\beta}(G)}$
with the corresponding one-sided approximation.

Let
\[
\mathcal{S}_1(G):=\{g({\mathbf x})={\mathbf a}{\mathbf x}^{\rm
t}+c\;:\;{\mathbf a}\in\RR^d,\, c\in\RR,\, {\mathbf x}\in G\}.
\]
The space $\mathcal{S}_1(G)$ will be the main approximation set in
this paper. Let also $\mathcal{T}=\{{\mathbf t}^1,\ldots,{\mathbf
t}^{d+1}\}$ be the $d$-dimensional simplex with vertices ${\mathbf
t}^j$, $j=1,\ldots,d+1$. We shall consider the following
optimization problem.

Let $Q({\mathbf x})=\textbf{xx}^{\rm t}$, and for
$\mathcal{T}\subset\RR^d$, set
\[
\sigma_{p;\alpha,\beta;d}(\mathcal{T}):=\frac{E(Q;\mathcal{S}_1(\mathcal{T}))_{L_{p;\alpha,\beta}(\mathcal{T})}}{|\mathcal{T}|^{1+\frac{1}{p}}},
\]
where $|\mathcal{T}|$ stands for the $d$-dimensional volume of the
simplex $\mathcal{T}$. The purpose of this paper is to solve
\begin{problem}\label{problem1}
Find
\begin{equation}\label{probl}
\sigma_{p;\alpha,\beta;d}:=\inf\limits_{\mathcal{T}}\sigma_{p;\alpha,\beta;d}(\mathcal{T})
\end{equation}
and describe all simplices $\mathcal{T}$, for which the infimum in
the right-hand part of~(2.3) is achieved.
\end{problem}

A solution to Problem~1 will allow to solve the following related problems.

\begin{problem}\label{problem2}
Find
\begin{equation}\label{probl'}
\sigma_{p;d}:=\inf\limits_{\mathcal{T}}\frac{E(Q;\mathcal{S}_1(\mathcal{T}))_{L_{p}(\mathcal{T})}}{|\mathcal{T}|^{1+\frac{1}{p}}}
\end{equation}
and describe all simplices $\mathcal{T}$, for which the infimum in
the right-hand part of~(2.4) is achieved.
\end{problem}

\begin{problem}\label{problem3}
Find
\begin{equation}\label{problI}
\sigma^{\pm}_{p;d}:=\inf\limits_{\mathcal{T}}\frac{E^{\pm}(Q;\mathcal{S}_1(\mathcal{T}))_{L_{p}(\mathcal{T})}}{|\mathcal{T}|^{1+\frac{1}{p}}}
\end{equation}
and describe all simplices $\mathcal{T}$, for which the infimum in
the right-hand part of~(\ref{problI}) is achieved.
\end{problem}

\section{History and the main result}\label{sec3}

Note that the quantity
$E^+(Q;\mathcal{S}_1(\mathcal{T}))_{L_{p}(\mathcal{T})}$ coincides
with the error of linear interpolation of the quadratic function
$Q({\mathbf x})$ at the vertices of the simplex $\mathcal{T}$, and
the quantity
$E^-(Q;\mathcal{S}_1(\mathcal{T}))_{L_{p}(\mathcal{T})}$ coincides
with the error of tangential interpolation of $Q({\mathbf x})$ on
the simplex $\mathcal{T}$. In view of formulas~(\ref{ter}), we
have
\[
\sigma^+_{p;d}=\lim\limits_{\beta\to+\infty}\sigma_{p;1,\beta;d}\quad
{\rm and}\quad
\sigma^-_{p;d}=\lim\limits_{\alpha\to+\infty}\sigma_{p;\alpha,1;d}.
\]

The quantity $\sigma^+_{p;d}$ has been considered for $d=2$ in
connection with the problem of finding the best triangulation
$\triangle_N$ consisting of $N$ triangles of the set
$G\subset\RR^2$, provided that the $L_p$-error of interpolation at
the vertices of $\triangle_N$ of convex function $f$ is minimized.

The first attempt to find $\sigma^+_{p;d}$ is due to D'Azevedo and
Simpson~\cite{Daz3}, who computed $\sigma_{\infty;2}^+$. To the
best of our knowledge, the progress on this problem can be
outlined as follows:
\begin{enumerate}
\item[(1)] $d=2$, $p=\infty$ \cite{Daz3};
\item[(2)] $d\ge 2$, $p=\infty$ \cite{Rajan_91};
\item[(3)] $d=3$, $p=2$ \cite{Brezin_92};
\item[(4)] $d=2$, $p=1$ \cite{boro};
\item[(5)] $d=2$, $p=2$ \cite{kodla1};
\item[(6)] $d\ge 2$, $p\in\NN$ \cite{chen_1};
\item[(7)] $d=2$, $p\in(0,\infty)$ \cite{BBS};
\item[(8)] $d\ge 2$, $p\in(1,\infty)$ \cite{chen}.
\end{enumerate}

\begin{remark}
Note that $\sigma^+_{p;d}$ for $d=2$ and $p\in(1,\infty)$ was
independently found by Chen~\cite{chen} and Babenko, {\it et al.}
\cite{BBS}.
\end{remark}

\begin{remark}
By $L_p$-error in the case $p\in(0,1)$, we
understand the following expression:
\[
E(f;H)_{L_p(G)}:=\inf\left\{\left(\int_G |f({\mathbf
x})-u({\mathbf x})|^p\,d{\mathbf x}\right)^{\frac{1}{p}}\;:\;u\in
H\right\}.
\]
\end{remark}

\begin{remark}
Infimum in Problem~\ref{problem3} is achieved only on regular
simplices.
\end{remark}

To the best of our knowledge, Problem~\ref{problem2} was solved
only in the case $d=2$, $p=2$ and $\alpha=\beta=1$ by
Nadler~\cite{Nadler_85,Nadler_86}. In the next section, we shall
give the solution of Problem~1 for all $\alpha,\beta>0$, $1\le
p\le\infty$ and $d\in\NN$.

Let $\mathcal{T}_0$ be a regular simplex of unit volume in
$\RR^d$. The main result of our paper is the following:
\setcounter{theorem}{0}
\begin{theorem}\label{theorem3.1}
Let $\alpha,\beta>0$, $d\in\NN$ and $1\le
p\le\infty$. Then
\[
\sigma_{p;\alpha,\beta;d}=\sigma_{p;\alpha,\beta;d}(\mathcal{T}_0).
\]
\end{theorem}
In view of~(\ref{ter}) we obtain the following statements.
\setcounter{theorem}{0}
\begin{corollary}\label{corollary3.1}
Let $d\in\NN$ and $1\le p\le\infty$.
Then
\[
\sigma_{p;d}=\frac{E(Q;\mathcal{S}_1(\mathcal{T}_0))_{L_p(\mathcal{T}_0)}}{|\mathcal{T}_0|^{1+\frac 1p}}.
\]
\end{corollary}

\begin{corollary}\label{corollary3.2}
Let $d\in\NN$ and $1\le p\le\infty$.
Then
\[
\sigma_{p;d}^{\pm}=\frac{E^{\pm}(Q;\mathcal{S}_1(\mathcal{T}_0))_{L_p(\mathcal{T}_0)}}{|\mathcal{T}_0|^{1+\frac 1p}}.
\]
\end{corollary}

Recall that $E^{+}(Q;\mathcal{S}_1(\mathcal{T}_0))$ is the error
of linear interpolation of\break \hbox{function}~$Q$.

The next section is devoted to the proof of this theorem.

\section{Proof of the main result}\label{sec4}

Let $\alpha,\beta>0$ be fixed throughout this section. Note that
the value of the quantity
$E(Q;\mathcal{S}_1(\mathcal{T}))_{L_{p;\alpha,\beta}(\mathcal{T})}$
is independent of translations of the simplex $\mathcal{T}$ and
its volume. For simplices $\mathcal{T},\mathcal{T}'\subset\RR^d$,
we shall write $\mathcal{T}=\mathcal{T}'$ if there exists a motion
$F$ of the space $\RR^d$ such that $F(\mathcal{T})=\mathcal{T}'$,
and we shall write $\mathcal{T}\ne \mathcal{T}'$ otherwise.

The proof of the main theorem consists of two parts contained in
the following two lemmas.

\begin{lemma}\label{lamme4.1}
Let $\mathcal{T}$ be an arbitrary
$d$-dimensional simplex of unit volume. Then there exists a
constant $C>0$, independent of $\mathcal{T}$, such that
\[
\sigma_{p;\alpha,\beta;d}(\mathcal{T})\ge C ({\rm diam}\,\mathcal{T})^2.
\]
\end{lemma}

\begin{lemma}\label{lamme4.2}
If $\mathcal{T}$, $\mathcal{T}\ne \mathcal{T}_0$, is a simplex of
unit volume in $\RR^d$ then there exists a simplex
$\mathcal{T}^*\subset\RR^d$ of unit volume such that
\[
\sigma_{p;\alpha,\beta;d}(\mathcal{T})>\sigma_{p;\alpha,\beta;d}(\mathcal{T}^*).
\]
\end{lemma}

Indeed, in view of Lemma~4.1, there exists an optimal
$d$-dimensional simplex $\mathcal{T}'$ of unit volume such that
$\sigma_{p;\alpha,\beta;d}=\sigma_{p;\alpha,\beta;d}(\mathcal{T}')$.
Then, Lemma~\ref{lamme4.2} gives $\mathcal{T}'=\mathcal{T}_0$.

\begin{proof}[Proof of Lemma 4.1]
Let $\mathcal{T}_d:=\mathcal{T}=\{{\mathbf t}^1,{\mathbf
t}^2,\ldots,{\mathbf t}^{d+1}\}$ be a simplex of unit volume.
Assume that $\|{\mathbf t}^1-{\mathbf t}^2\|_2={\rm
diam}\,\mathcal{T}_d$. In addition, for $j=1,\ldots,d-1$, let
$\mathcal{T}_{j}=\{{\mathbf t}^1,{\mathbf t}^2,\ldots, {\mathbf
t}^{j+1}\}$ be a simplex in $\RR^{j}$.

First, let us consider the case $1\le p<\infty$. For
$j=2,\ldots,d$, by $h_j$ denoting the length of the height from
the vertex ${\mathbf t}^j$ of the simplex $\mathcal{T}_j$ to the
simplex $\mathcal{T}_{j-1}$. For every ${\mathbf
a}=(a_1,\ldots,a_d)\in\RR^d$, let ${\mathbf
a}'=(a_1,\ldots,a_{d-1})$. Then,
\begin{align*}
\sigma_{p;\alpha,\beta;d}^p(\mathcal{T}_d) & =  \inf_{{\mathbf
a}\in\RR^d,\; c\in\RR}
\int_{\mathcal{T}_d}|\textbf{xx}^{\rm t}-{\mathbf a}{\mathbf x}^{\rm t}-c|^p_{\alpha,\beta}\,d{\mathbf x}\\
& =  \inf_{{\mathbf a}\in\RR^d,\; c\in\RR}\int_{0}^{h_d}\int_{\frac{u}{h_d}\mathcal{T}_{d-1}}|u^2+{\mathbf y}{\mathbf y}^{\rm t}-a_d u-{\mathbf a}'{\mathbf y}^{\rm t}-c|^p_{\alpha,\beta}\,d{\mathbf y}du \\
&\ge
\int_{0}^{h_d}\left(\frac{u}{h_d}\right)^{2p+(d-1)}\sigma_{p;\alpha,\beta;d-1}^p(\mathcal{T}_{d-1})\,du\\
&=\frac{h_d}{2p+d}\,\sigma_{p;\alpha,\beta;d-1}^p(\mathcal{T}_{d-1}).
\end{align*}
Proceeding by induction on $d$, we verify that
\begin{align*}
\sigma_{p;\alpha,\beta;d}^p(\mathcal{T}_{d}) & \ge  \frac{h_2
h_3\ldots h_d}{(2p+2)(2p+3)\cdots(2p+d)}
\sigma^p_{p;\alpha,\beta;1}(\mathcal{T}_1) \\ &\ge  C^p\cdot
\frac{h_2 h_3\ldots h_d }{d!}\cdot({\rm
diam}\,\mathcal{T}_{d})^{2p+1} =C^p ({\rm
diam}\,\mathcal{T}_{d})^{2p},
\end{align*}
where $\Upsilon$ is some positive constant independent of the simplex $\mathcal{T}_{d}$.

Let us turn to the case $p=\infty$. In this case, we obtain
\begin{align*}
\sigma_{\infty;\alpha,\beta;d}(\mathcal{T}_d) & \ge
\sigma_{\infty;\alpha,\beta;1}(\mathcal{T}_1) \\ & =
\inf_{k\in\RR,c\in\RR}\sup\limits_{u\in[0,{\rm
diam}\,\mathcal{T}_{d}]}|u^2-ku-c|_{\alpha,\beta}\ge\Upsilon({\rm
diam}\,\mathcal{T}_d)^2.
\end{align*}
\vspace*{-1pc}\end{proof}

\pagebreak
\begin{proof}[Proof of Lemma 4.2]
Let $\mathcal{T}=\{{\mathbf w}^1,{\mathbf w}^2,{\mathbf
t}^1,\ldots, {\mathbf t}^{d-1}\}\ne \mathcal{T}_0$ be a simplex of
unit volume. Without loss of generality, we may assume that
$\|{\mathbf w}^1-{\mathbf t}^1\|_2\ne\|{\mathbf w}^2-{\mathbf
t}^1\|_2$. Clearly, we can always choose the coordinate system in
$\RR^d$ so that the vertices of $\mathcal{T}$ have the following
coordinates:
\[
{\mathbf w}^1=(-\delta,0,0,\ldots,0),\quad {\mathbf w}^2=(\delta,0,0,\ldots,0),
\]
where $\delta:=\frac{1}{2}\|{\mathbf w}^1-{\mathbf w}^2\|_2$, and
the remaining vectors ${\mathbf t}^1,{\mathbf t}^2,\ldots,{\mathbf
t}^{d-1}$ have the following coordinates:
\[
\begin{array}{@{}lcccccccccl@{}}
{\mathbf t}^1 & =: & ( & b_1 & a_{1,1} & 0 & 0 & \ldots & 0 & 0 & ), \\
{\mathbf t}^2 & =: & ( & b_2 & a_{1,2} & a_{2,2} & 0 & \ldots & 0 & 0 & ), \\
{\mathbf t}^3 & =: & ( & b_3 & a_{1,3} & a_{2,3} & a_{3,3} & \ldots & 0 & 0 & ), \\
\vdots & & & \vdots & \vdots & \vdots & \vdots & & \vdots & \vdots & \\
{\mathbf t}^{d-2} & =: & ( & b_{d-2} & a_{1,d-2} & a_{2,d-2} & a_{3,d-2} & \ldots & a_{d-2,d-2} & 0 & ), \\
{\mathbf t}^{d-1} & =: & ( & b_{d-1} & a_{1,d-1} & a_{2,d-1} & a_{3,d-1} & \ldots & a_{d-2,d-1} & a_{d-1,d-1} & ).
\end{array}
\]
Note that in view of our assumption, $b_1\ne 0$. In addition, it can be easily seen that $a_{j,j}\ne 0$ for all $j=1,\ldots,d-1$.

Let
\[
{\mathbf b}:=(b_1,\ldots,b_{d-1})
\]
and
\[
{\mathbf A}:=\left(\begin{array}{@{}cccccc@{}}
a_{1,1} & a_{1,2} & a_{1,3} & \ldots & a_{1,d-2} & a_{1,d-1} \\
0 & a_{2,2} & a_{2,3} & \ldots & a_{2,d-2} & a_{2,d-1} \\
0 & 0 & a_{3,3} & \ldots & a_{3,d-2} & a_{3,d-1} \\
\vdots& \vdots & \vdots & \ddots & \vdots & \vdots \\
0 & 0 & 0 & \ldots &a_{d-2,d-2}& a_{d-2,d-1} \\
0 & 0 & 0 & \ldots & 0 & a_{d-1,d-1}
\end{array}\right).
\]
Since the matrix ${\mathbf A}$ is non-singular, set
\begin{equation}\label{dfd}
{\mathbf y}=(y_1,\ldots,y_{d-1}):={\mathbf b A}^{-1}.
\end{equation}

Let ${\mathbf I}$ be the identity matrix of size $(d-1)\times
(d-1)$. In addition, set ${\mathbf R}:={\mathbf y}^{\rm t}
{\mathbf y}+{\mathbf I}$. It can be easily seen that the matrix
${\mathbf R}$ is positive definite. Therefore, in view of the
Cholesky decomposition (a standard technique in numerical
analysis, whose description can be found, for instance,
in~\cite{Num_An}), there exists an upper triangular matrix
${\mathbf U}=(u_{k,j})_{1\le k,j\le d-1}$ such that
\[
{\mathbf R = \textbf{U}}^{\rm t}{\mathbf U}.
\]
Moreover, for every $j=1,\ldots,d-1$, we have $u_{j,j}=\sqrt{\frac{D_{j}}{D_{j-1}}}$, where
\begin{equation}\label{zzz}
D_0:=1,\quad D_k:=\det\left(\begin{array}{@{}cccccc@{}}
1+y_1^2 & y_1y_2 & y_1y_3 & \ldots & y_1y_{k-1} & y_1y_{k} \\
y_1y_2 & 1+y_2^2 & y_2y_3 & \ldots & y_2y_{k-1} & y_2y_{k} \\
y_1y_3 & y_2y_3 & 1+y_3^2 & \ldots & y_3y_{k-1} & y_3y_{k} \\
\vdots & \vdots & \vdots & \ddots & \vdots & \vdots \\
y_1y_{k-1} & y_2y_{k-1} & y_3y_{k-1} & \ldots & 1+y_{k-1}^2 & y_{k-1}y_{k} \\
y_1y_{k} & y_2y_{k} & y_3y_{k} & \ldots & y_{k-1}y_{k} & 1+y_{k}^2
\end{array}\right),
\end{equation}
for all $ k=1,\ldots,d-1$. Consequently, $u_{j,j}>0$ for every
$j=1,\ldots,d-1$.

Let ${\mathbf Q}=(q_{k,j})_{1\le k,j\le d-1}$ be the diagonal
matrix such that $q_{j,j}=u_{j,j}$, $j=1,\ldots,d-1$. We define
\[
\overline{\mathbf U}:={\mathbf Q}^{-1}{\mathbf U}.
\]
Then $\overline{\mathbf U}$ is the unit upper triangular matrix.
Therefore,
\[
{\mathbf R}=\overline{\mathbf U}^{\rm t}{\mathbf Q}^2\overline{\mathbf U}.
\]
Set
\[
{\mathbf M}:=\overline{\mathbf U}{\mathbf A}.
\]
Denote the elements of ${\mathbf M}$ by $m_{k,j}$, i.e., ${\mathbf
M}=(m_{k,j})_{1\le k,j\le d-1}$. Note that for every
$j=1,\ldots,d-1$, it follows that $m_{j,j}=a_{j,j}$.

Let us now consider the simplex
$\widetilde{\mathcal{T}}=\{{\mathbf w}^1,{\mathbf
w}^2,\widetilde{\mathbf t}^1, \ldots, \widetilde{\mathbf
t}^{d-1}\}$, whose vertices have the following coordinates:
\[
\begin{array}{@{}lcccccccccl@{}}
\widetilde{\mathbf t}^1 & := & ( & 0 & a_{1,1} & 0 & 0 & \ldots & 0 & 0 & ), \\
\widetilde{\mathbf t}^2 & := & ( & 0 & m_{1,2} & a_{2,2} & 0 & \ldots & 0 & 0 & ), \\
\widetilde{\mathbf t}^3 & := & ( & 0 & m_{1,3} & m_{2,3} & a_{3,3} & \ldots & 0 & 0 & ), \\
\vdots & & & \vdots & \vdots & \vdots & \vdots & \vdots & & \vdots & \\
\widetilde{\mathbf t}^{d-2} & := & ( & 0 & m_{1,d-2} & m_{2,d-2} & m_{3,d-2} & \ldots & a_{d-2,d-2} & 0 & ), \\
\widetilde{\mathbf t}^{d-1} & := & ( & 0 & m_{1,d-1} & m_{2,d-1} & m_{3,d-1} & \ldots & m_{d-2,d-1} & a_{d-1,d-1} & ).
\end{array}
\]
Obviously, the volumes of simplices $\mathcal{T}$ and
$\widetilde{\mathcal{T}}$ coincide.

Now we shall construct the linear transformation $S:\RR^d\to\RR^d$
such that $S(\widetilde{\mathcal{T}})=\mathcal{T}$. To that end,
set
\[
{\mathbf h}:={\mathbf b}{\mathbf M}^{-1}\in\RR^{d-1},
\]
and define the linear transformation $S$ with the help of the matrix
\[
{\mathbf S}=\left(\begin{matrix} 1 & {\mathbf h} \\ \begin{matrix} 0
\\ \vdots \\ 0 \end{matrix} & (\overline{\mathbf U})^{-1}\end{matrix}\right).
\]
Note that $\det{\mathbf S}=1$. Then for every ${\mathbf
a}=(a_1,a_2,\ldots,a_d)\in\RR^d$ and $c\in\RR$, it follows that
\begin{equation}\label{equa1}
L:=\|\textbf{xx}^{\rm t}-{\mathbf a}{\mathbf x}^{\rm
t}-c\|_{L_{p;\alpha,\beta}(\mathcal{T})}=\|{\mathbf v}{\mathbf
S}^{\rm t}{\mathbf S}{\mathbf v}^{\rm t}-{\mathbf a}{\mathbf
S}{\mathbf v}^{\rm
t}-c\|_{L_{p;\alpha,\beta}(\widetilde{\mathcal{T}})}.
\end{equation}

Let us consider the simplex $\widehat{\mathcal{T}}=\{{\mathbf
w}^1,{\mathbf w}^2,\widehat{\mathbf t}^1,\ldots, \widehat{\mathbf
t}^{d-1}\}$ such that
\[
\begin{array}{@{}lcccccccccl@{}}
\widehat{\mathbf t}^1 & := & ( & -b_1 & a_{1,1} & 0 & 0 & \ldots & 0 & 0 & ), \\
\widehat{\mathbf t}^2 & := & ( & -b_2 & a_{1,2} & a_{2,2} & 0 & \ldots & 0 & 0 & ), \\
\widehat{\mathbf t}^3 & := & ( & -b_3 & a_{1,3} & a_{2,3} & a_{3,3} & \ldots & 0 & 0 & ), \\
\vdots & & & \vdots & \vdots & \vdots & \vdots & & \vdots & \vdots & \\
\widehat{\mathbf t}^{d-2} & := & ( & -b_{d-2} & a_{1,d-2} & a_{2,d-2} & a_{3,d-2} & \ldots & a_{d-2,d-2} & 0 & ), \\
\widehat{\mathbf t}^{d-1} & := & ( & -b_{d-1} & a_{1,d-1} & a_{2,d-1} & a_{3,d-1} & \ldots & a_{d-2,d-1} & a_{d-1,d-1} & ).
\end{array}
\]
It can be easily verified that the linear transformation $\widehat{S}:\RR^d\to\RR^d$ defined with the help of the matrix
\[
\widehat{\mathbf S}=\left(\begin{matrix} 1 & -{\mathbf h} \\
\begin{matrix} 0 \\ \vdots \\ 0 \end{matrix} & (\overline{\mathbf
U})^{-1}\end{matrix}\right)
\]
transforms the simplex $\widetilde{\mathcal{T}}$ into the simplex
$\widehat{\mathcal{T}}$. Therefore, for $\widehat{\mathbf
a}:=(-a_1,$\break $a_2,\ldots,a_d)$, we obtain
\begin{equation}\label{equa2}
\widehat{L}:=\|\mathbf{xx}^{\rm t}-\widehat{\mathbf a}{\mathbf
x}^{\rm
t}-c\|_{L_{p;\alpha,\beta}(\widehat{\mathcal{T}})}=\|{\mathbf
v}\widehat{\mathbf S}^{\rm t}\widehat{\mathbf S}{\mathbf v}^{\rm
t}-\widehat{\mathbf a}\widehat{\mathbf S}{\mathbf v}^{\rm
t}-c\|_{L_{p;\alpha,\beta}(\widetilde{\mathcal{T}})}.
\end{equation}

Due to the symmetry of simplices $\mathcal{T}$ and
$\widehat{\mathcal{T}}$, we obtain $L=\widehat{L}$. Then
from~(\ref{equa1}) and~(\ref{equa2}), we derive that
\begin{equation}\label{main_equa}
L=\frac{1}{2}(L+\widehat{L})\ge\left\|\frac{1}{2}\left[{\mathbf
v}({\mathbf S}^{\rm t}{\mathbf S}+\widehat{\mathbf S}^{\rm
t}\widehat{\mathbf S}){\mathbf v}^{\rm t}-({\mathbf a}{\mathbf
S}+\widehat{\mathbf a}\widehat{\mathbf S}){\mathbf v}^{\rm
t}-2c\right]\right\|_{L_{p;\alpha,\beta}(\widetilde{\mathcal{T}})}.
\end{equation}

Note that
\[
\begin{array}{@{}rcl@{}}
2{\mathbf D}:={\mathbf S}^{\rm t}{\mathbf S}+\widehat{\mathbf S}^{\rm t}\widehat{\mathbf S} & = &
\left(\begin{array}{@{}cc@{}} 1 & \begin{array}{@{}ccc@{}} 0 & \ldots & 0 \end{array} \\ {\mathbf h}^{\rm t}
&\left[(\overline{\mathbf U})^{-1}\right]^{\rm t}\end{array}\right)\left(\begin{array}{@{}cc@{}}
1 & {\mathbf h} \\ \begin{array}{c} 0 \\ \vdots \\ 0 \end{array}
& (\overline{\mathbf U})^{-1}\end{array}\right) \\
&&+\left(\begin{array}{@{}cc@{}} 1 & \begin{array}{ccc} 0 & \ldots & 0
\end{array} \\ -{\mathbf h}^{\rm t} &
\left[(\overline{\mathbf U})^{-1}\right]^{\rm t}\end{array}\right)
\left(\begin{array}{@{}cc@{}} 1 & -{\mathbf h} \\
\begin{array}{@{}c@{}} 0 \\ \vdots \\ 0 \end{array} & (\overline{\mathbf U})^{-1}\end{array}\right) \\
&= & 2\left(\begin{array}{@{}cc@{}}1 & \begin{array}{@{}ccc@{}} 0 & \ldots & 0
\end{array} \\ \begin{array}{@{}c@{}} 0 \\ \vdots \\ 0 \end{array} &
{\mathbf h}^{\rm t}{\mathbf h}+[(\overline{\mathbf
U})^{-1}]^{\rm t}(\overline{\mathbf U})^{-1}
\end{array}\right).
\end{array}
\]
Since ${\mathbf h}={\mathbf b}{\mathbf M}^{-1}={\mathbf b}{\mathbf
A}^{-1}(\overline{\mathbf U})^{-1}$, we have
\begin{align*}
{\mathbf h}^{\rm t}{\mathbf h}+[(\overline{\mathbf
U})^{-1}]^{\rm t}(\overline{\mathbf U})^{-1} & =
[(\overline{\mathbf U})^{-1}]^{\rm t}[({\mathbf
A}^{-1})^{\rm t}{\mathbf b}^{\rm t}{\mathbf b}{\mathbf
A}^{-1}+{\mathbf I}](\overline{\mathbf U})^{-1} \\
&=[(\overline{\mathbf U})^{-1}]^{\rm t}{\mathbf
R}(\overline{\mathbf U})^{-1}={\mathbf Q}^2.
\end{align*}
Therefore,
\[
{\mathbf D}=\left(\begin{array}{@{}cc@{}}1 & \begin{array}{@{}ccc@{}} 0 & \ldots & 0 \end{array} \\
\begin{array}{@{}c@{}} 0 \\ \vdots \\ 0 \end{array} & {\mathbf Q}^2 \end{array}\right),
\]
and ${\mathbf D}$ is the diagonal matrix with elements
$1,D_1,\frac{D_2}{D_1},\ldots,\frac{D_{d-1}}{D_{d-2}}$ on the main
diagonal (numbers $D_j$, $j=1,\ldots,d-1$, were defined
in~(\ref{zzz})). Let $F$ be the linear transformation defined with
the help of the diagonal matrix ${\mathbf F}$, having elements
$D_{d-1}^{\frac 1{2d}},D_{d-1}^{\frac 1{2d}}\cdot
\sqrt{\frac{1}{D_1}},D_{d-1}^{\frac 1{2d}}\cdot
\sqrt{\frac{D_1}{D_2}},\ldots,D_{d-1}^{\frac
1{2d}}\cdot\sqrt{\frac{D_{d-2}}{D_{d-1}}}$ on the main diagonal.
It can be easily seen that $\det{\mathbf F}=1$. Let
$\mathcal{T}^*=F^{-1}(\widetilde{\mathcal{T}})$. Then, in view
of~(\ref{main_equa}),
\[
L\ge D_{d-1}^{\frac 1d}\|{\mathbf z}{\mathbf z}^{\rm t}-{\mathbf g}{\mathbf z}^{\rm t}-c'\|_{L_{p;\alpha,\beta}(\mathcal{T}^*)},
\]
where ${\mathbf g}=\frac{1}{2D_{d-1}^{\frac 1{2d}}}({\mathbf
a}{\mathbf S}+\widehat{\mathbf a}\widehat{\mathbf S}){\mathbf F}$
and $c'=\frac{c}{D_{d-1}^{\frac 1{2d}}}$. Consequently,
\[
\sigma_{p;\alpha,\beta;d}(\mathcal{T})\ge D_{d-1}^{\frac
1d}\sigma_{p;\alpha,\beta;d}(\mathcal{T}^*).
\]

Let us show that the assumption $\|{\mathbf w}^1-{\mathbf
t}^1\|_2\ne\|{\mathbf w}^2-{\mathbf t}^1\|_2$ yeilds\break $D_{d-1}>1$.
Indeed, since the matrix ${\mathbf R}$ is positive definite, it
follows that
\begin{equation}\label{trtt}
D_{d-1}=\det{({\mathbf y}^{\rm t}{\mathbf
y}+I)}=(1+\lambda_1)\cdots(1+\lambda_{d-1})\ge 1,
\end{equation}
where $\lambda_j\ge 0$, $j=1,\ldots,d-1$, are the eigenvalues of
the matrix ${\mathbf y}^{\rm t}{\mathbf y}$ (the existence of
non-negative eigenvalues is guaranteed by positive
semidefiniteness and symmetry of ${\mathbf y}^{\rm t}{\mathbf
y}$). Inequality~(\ref{trtt}) becomes an equality if and only if
$\lambda_j=0$ for all $j=1,\ldots,d-1$, or, that is the same,
$y_j=0$ for every $j=1,\ldots,d-1$. In view of~(\ref{dfd}) and
non-singularity of the matrix ${\mathbf A}$, this yields that
$b_j=0$ for all $j=1,\ldots,d-1$. However, this contradicts the
fact that $b_1\ne 0$, which follows from the assumption
$\|{\mathbf w}^1-{\mathbf t}^1\|_2\ne\|{\mathbf w}^2-{\mathbf
t}^1\|_2$.
\end{proof}

\section*{Acknowledgments}

The results have been obtained during the visit of Babenko,
Parfinovych and Skorokhodov to Sam Houston State University. V. Babenko was supported by SHSU Enhancement Grant for Research.

%\addres{Ul. Kazakova, d. 4A, kv. 33\\
%Dnepropetrovsk 49050\\
%Ukraine\\
%\email{\textit{E-mail address}: {\tt babenko@shsu.edu}}\\
%\received{April 02, 2009}}


\begin{thebibliography}{99}
\bibitem{Bab_82} V.F. Babenko, \textit{Non-symmetric approximations in spaces of summable functions}, Ukrain. Mat. Zh. \textbf{34} (1982), 409--416; English transl.: Ukrainian Math. J. \textbf{34} (1982), 331--336.

\bibitem{Bab_83} V.F. Babenko, \textit{Asymmetric extremal problems in approximation theory}, Dokl. SSSR \textbf{269}(3) (1983), 521--524 (in Russian).

\bibitem{Bab_84} V.F. Babenko, \textit{Duality theorems for some problems in approximation theory}, Contemp. questions of real and complex analysis, In-tyt math. AN USSR, Kiev (1984), 3--13 (in Russian).

\bibitem{Bab_87} V.F. Babenko, \textit{Approximations, widths and optimal quadrature formulae for classes of periodic functions with rearrangement invariant sets of derivatives}, Anal. Math. \textbf{13} (1987), 15--28.

\bibitem{BBLS} V. Babenko, Yu. Babenko, A. Ligun and A. Shumeiko, \textit{On asymptotical behavior of the optimal linear spline interpolation error of} $C^2$ \textit{functions}, East J. Approx. \textbf{12}(1) (2006), 71--101.

\bibitem{BBS} V. Babenko, Yu. Babenko and D. Skorokhodov, \textit{Exact asymptotics of the optimal} $L_{p,\Omega}$\textit{-error of linear spline interpolation}, East J. Approx. \textbf{14}(3) (2008), 285--317.

\bibitem{boro} K. B\"{o}r\"{o}czky and M. Ludwig, \textit{Approximation of convex bodies and a momentum lemma for power diagrams}, Monatsh. Math., \textbf{127}(2) (1999),\break 101--110.

\bibitem{Boro1} K. B\"{o}r\"{o}czky, \textit{Approximation of general smooth convex bodies}, Adv. Math., \textbf{153} (2000), 325--341.

\bibitem{Brezin_92} M. Brezin, \textit{A solution-based triangular and thetrahedral mesh quality indicator}, SIAM J. Sci. Comput. \textbf{19} (1992), 979--997.

\bibitem{chen_1} L. Chen, P. Sun and J. Xu, \textit{Optimal anisotropic meshes for minimizing interpolation errors in} $L_p$-\textit{norm}, Math. Comp. \textbf{76} (2007), 179--204.

\bibitem{chen} L. Chen, \textit{On minimizing the linear interpolation error of convex quadratic functions and the optimal simplex}, East J. Approx. \textbf{14}(3) (2008), 271--284.

\bibitem{Daz3} E.F. D'Azevedo and R.B. Simpson, \textit{On optimal interpolation triangle incidences}, SIAM J. Sci. Stat. Comput. \textbf{10}(6) (1989), 1063--1075.

\bibitem{Nira_Dyn_1} N. Dyn, D. Revin and S. Rippa, \textit{Data dependent triangulations for piecewise linear interpolation}, IMA, J. Numer. Anal. \textbf{10} (1988), 137--154.

\bibitem{Nira_Dyn_2} N. Dyn, D. Revin and S. Rippa, \textit{Algorithms for construction of data dependent triangulations}, in `Algorithms for approximation II', eds. J.C. Mason and M.G. Cox, Chapman and Hall, New York, 1990, 185--192.

\bibitem{Num_An} C. Gerald and P. Wheatley, Applied Numerical Analysis, 7th edn., Addison-Wesley, Reading, MA, 2003.

\enlargethispage{20pt}\bibitem{Gruber} P. Gruber, \textit{Error of asymptotic formulae for volume approximation of convex bodies in} $E^d$, Monatsh. Math. \textbf{135} (2002), 279--304.

\bibitem{Handbook} J.E. Goodman and J.O'Rourke, \textit{Handbook of discrete and computational geometry}, in Discrete mathematics and its applications, 2nd edn., CRC Press LLC, Boca Raton, FL, 2004.

\bibitem{huang} W. Huang, \textit{Measuring mesh qualities and application to variational mesh adaption}, SIAM J. Sci. Comput. \textbf{26}(5) (2005), 1643--1666.

\bibitem{KLD} N.P. Korneichuk, A.A. Ligun and V.G. Doronin, \textit{Approximation with constraints}, Naukova dumka, Kiev, 1982 (in Russian).

\bibitem{Korn} N.P. Korneichuk, \textit{Exact constants in approximation theory}, Nauka, Moscow, 1987; translated from the Russian by K. Ivanov. \textit{Encyclopedia of mathematics and its applications}, \textbf{38}, Cambridge University Press, Cambridge, 1991.

\bibitem{Krein} M.G. Krein, \textit{The L-problem in an abstract linear normed space}, in Some questions in the theory of moments, eds. N.I. Akhiezer and M.G. Krein, American Mathematical Society, Providence, RI, 1962, 175--204.

\bibitem{KreinNudelman} M.G. Krein and A.A. Nudel'man, \textit{The Markov moment problem and extremal problems}, Nauka, Moscow, 1973; English transl.: Fundamentals of the theory of quasigroups and loops, American Mathematical Society, Providence, 1977.

\bibitem{Nadler_85} E. Nadler, \textit{Piecewise linear approximation on triangulations of a planar region}, PhD thesis, Brown University, 1985.

\bibitem{Nadler_86} E. Nadler, \textit{Piecewise linear best} $L_2$ \textit{approximation on triangulations}, in Approximation theory, eds. C.K. Chui, L.L. Schumaker and J.D. Ward, \textbf{V}, 1986, 499--502.

\bibitem{kodla1} H. Pottmann, R. Krasauskas, B. Hamann, K. Joy and W. Seibold, \textit{On piecewise linear approximation of quadratic functions}, J. Geom. Graph. \textbf{4}(1) (2000), 31--53.

\bibitem{Rajan_91} V.T. Rajan, \textit{Optimality of Delaunay triangulation in} $\RR^d$, Proceedings of the Seventh Annual Symposium on Comp. Geom. (1991), 357--363.

\bibitem{Toth} L. Fejes Toth, \textit{Lagerungen in der Ebene, auf der Kugel und im Raum}, 2nd edn. Springer, Berlin, 1972.
\end{thebibliography}
\end{document}